\documentclass[12pt]{article}

\usepackage{amsfonts}

\title{Square Eulerian Quadruples}
\author{Allan J. MacLeod\\Dept. of Mathematics and Statistics\\
University of Paisley\\High St.\\Paisley\\Scotland\\PA1 2BE\\
(e-mail: allan.macleod@paisley.ac.uk)}
\date{}

\begin{document}

\maketitle

\begin{abstract}
We consider the problem of finding four different rational
squares, such that the product of any two plus the sum of the same
two always gives a square. We give some historical background to
this problem and exhibit one such quadruple.
\end{abstract}

\section{Introduction}
The genesis of this work is the following section taken from
Chapter $19$ in Volume $2$ of Dickson's "History of the Theory of
Numbers" \cite{dick}.

\bigskip

{\it Fermat treated the problem to find four numbers such that
the product of any two increased by the sum of those two gives a
square. He made use of three squares such that the product of any
two increased by the sum of the same two gives a square. Stating
that there is an infinitude of such sets of three squares, he
cited $4, 3504384/d, 2019241/d$, where $d=203401$. However, he
actually used the squares $25/9, 64/9, 196/9$, of Diophantus V,
5, which have the additional property that the product of any two
increased by the third gives a square. Taking these three squares
as three of our numbers and $x$ as the fourth, we are to satisfy
\begin{displaymath}
\frac{34}{9}x+\frac{25}{9}=\square \hspace{2cm}
\frac{73}{9}x+\frac{64}{9}=\square \hspace{2cm}
\frac{205}{9}x+\frac{196}{9}=\square
\end{displaymath}
This "triple equation" with squares as constant terms is readily
solved. T.L. Heath found $x$ to be the ratio of two numbers each
of 21 digits. }

\bigskip
This section generated several questions:
\begin{enumerate}
\item Is there "an infinitude of such sets of three squares"?
\item Where did the example $4, 3504384/d, 2019241/d$ come from?
\item What is the $x$ found by Heath and can we find a smaller
value - smaller meaning fewer digits in its rational form?
\item What would the solution have been with the original example
of Fermat?
\item Can we use the set of $3$ squares to find a fourth square,
so that the $4$ squares are a solution to Fermat's problem?
\end{enumerate}

Given a set $\{a_1,a_2,\ldots,a_m\}$ such that $a_i
a_j+a_i+a_j=\square$ for $1 \le i < j \le m$, Dujella \cite{duj1}
\cite{duj2} called such a set an {\bf Eulerian m-tuple}, in honour
of Euler who found the first example of a quadruple of numbers.

To be historically accurate, however, we should note that what
Euler was considering was sets $\{b_1,b_2,b_3,b_4\}$ with $b_i
b_j-1=\square$. This is related to the present problem since
\begin{displaymath}
a_i a_j + a_i + a_j = ( a_i + 1 )( a_j + 1) -1
\end{displaymath}

\section{Infinitude of Eulerian Triples of Squares}
We wish to find $x,y,z$ such that
\begin{eqnarray*}
x^2 y^2 + x^2 + y^2 & = & a^2 \\
x^2 z^2 + x^2 + z^2 & = & b^2 \\
y^2 z^2 + y^2 + z^2 & = & c^2 \\
\end{eqnarray*}
Rather than have a complicated analysis, we show that there are an
infinite number of solutions with $x=2$.

Thus we look for $y,z$ with
\begin{equation}
5y^2+4 =a^2 \hspace{1cm} 5z^2+4=b^2 \hspace{1cm} y^2
z^2+y^2+z^2=c^2
\end{equation}

Rational solutions of $5t^2+4=\square$ can be parameterised as
$t=4f/(5-f^2)$, with $f$ rational. Thus $y=4m/(5-m^2)$ and
$z=4n/(5-n^2)$, which can be substituted into the third part of
(1).

We find that $m$ and $n$ must satisfy
\begin{displaymath}
d^2=m^2 n^4+(m^4-4m^2+25)n^2+25m^2
\end{displaymath}
for $d$ rational.

Define $d=Y/m$ and $n=X/m$, giving the quartic relation
\begin{displaymath}
Y^2= X^4+(m^4-4m^2+25)X^2+25m^4
\end{displaymath}

This quartic has an obvious rational point $X=0, Y=5m^2$, and so
is birationally equivalent to an elliptic curve. Using the
standard method described by Mordell \cite{mord}, we find the
elliptic curve (with $m=p/q$)
\begin{equation}
J^2=K^3-2(p^4-4p^2q^2+25q^4)K^2+(p^8-8p^6q^2-34p^4q^4-200p^2q^6+625q^8)K
\end{equation}
with the relation $n=J/(2 p q K)$.

The curve has $3$ points of order $2$ with $J=0$, which lead to
$z=0$, so we need other rational points for non-trivial
solutions. Evidence suggests that the torsion subgroup is just
isomorphic to $\mathbb{Z}2 \times \mathbb{Z}2$, but this would be
difficult to prove. If true, we need the curves to have rank
greater than $0$ for solutions.

If we experiment, we quickly find that $m=1/2$ gives a curve
$J^2=K^3-770K^2+146625K$ which has rank $1$ with generator
$P=(245,2100)$. $m=1/2$ gives $y=8/19$ and the point P gives
$n=15/7$ and $z=21$, and it easily checked that
$\{4,64/361,441\}$ is a square Eulerian triple.

Since the curve has rank $1$ we have that integral multiples of P
are also rational points on the curve. For example, doubling the
point P leads to the values $K=187489/441, J=651232/9261$ and so
$n=376/9093$ and an alternative $z$ of $13675872/413271869$.
There are thus an infinite number of Eulerian triples with
$x=2,y=8/19$.

\section{Numerical Values}
It is totally unclear from Fermat's original work where his
numerical example involving $4$ comes from, as he just states
this result with no supporting algebra or computation. As we saw
in this last section, there are much simpler sets which include
$4$.

We now consider the system
\begin{displaymath}
\frac{34}{9}x+\frac{25}{9}=\square \hspace{1.5cm}
\frac{73}{9}x+\frac{64}{9}=\square \hspace{1.5cm}
\frac{205}{9}x+\frac{196}{9}=\square
\end{displaymath}

We first write
\begin{displaymath}
\frac{34}{9}x+\frac{25}{9}=\left( \frac{5}{3}+f x \right)^2
\end{displaymath}
which gives $x=2(17-15f)/(9f^2)$.

Substituting this into the second and third equations, we find
the following two equations must have rational solutions:
\begin{eqnarray*}
2(288f^2-1095f+1241)&=e^2 \\
2(882f^2-3075f+3485)&=h^2
\end{eqnarray*}

Now, in the first equation, $f=17/15$ gives $e=136/5$, so, if we
substitute $e=136/5+g(f-17/15)$, we can solve to find
\begin{displaymath}
f=\frac{17g^2-816g-23058}{15(g^2-576)}
\end{displaymath}
and if we substitute this into the other quadratic relation we
find, clearing denominators, that we must have a rational
solution to
\begin{equation}
s^2=14161g^4+731544g^3+28206441g^2+639486144g+6471280836=G(g)
\end{equation}

Since $14161=119^2$ we can attempt to complete the square, by
forming $G(g)-(\alpha g^2 +\beta g+\gamma)^2$. Simple arithmetic
shows that if $\alpha=119, \beta=21516/7, \gamma=919177353/11662$
then
\begin{displaymath}
G(g)-(\alpha g^2+\beta g+\gamma)^2 =
\frac{1581221503125(13328g+22275)}{136002244}
\end{displaymath}

Thus, we have a square solution when $g=-22275/13328$, which
gives $f=142415972261/56567733755$ and finally
\begin{displaymath}
x=\frac{-459818598496844787200}{631629004828419699201}
\end{displaymath}
where both numerator and denominator have $21$ digits.

Because equation (3) has leading coefficient $119^2$ it can be
transformed into an elliptic curve. Mordell's method leads to the
curve $j^2=k^3+20478k^2+99801585k$ with
\begin{displaymath}
g=\frac{3(105j-44(163k+2348685)}{34(49k+560880)}
\end{displaymath}

This curve has $3$ points of order $2$ with $k=0, -7999, -12483$
all with $j=0$. The first leads to $x=0$ whilst the other two
give undefined values for $f$. Thus we need points of infinite
order. Cremona's {\bf mwrank} program shows the rank to be $2$
with generators $(-9984,-222768)$ and $(-8379,-114912)$. The
first gives $g=-543/8$, $f=269/147$ and $x=-50176/72361$,
significantly simpler than Heath's value.

With regard to the $\{4, 3504384/203401, 2019241/203401\}$ triple
of squares we can perform an identical analysis. The
corresponding value of $x$ from completing the square is
\begin{displaymath}
x=\frac{-28448417598272924003671204878289354665765410185967616}
{36828906078832095599985737816846193226885934523284161}
\end{displaymath}
where both numerator and denominator have $53$ digits.

Attempts, as before, to find a smaller value of $x$ lead to the
elliptic curve
\begin{displaymath}
v^2=u^3+10450883424805u^2+26734915668323655104674200u
\end{displaymath}

The completing the square values lead to a point of infinite order
on this curve with $u=-9390695817653070336/2019241$.
Investigations with APECS, mwrank, and SAGE were unable to find
other generators. Both APECS and mwrank give $3$ as an upper
bound for the rank. The root number is $-1$, so the parity
conjecture suggests the rank is $1$ or $3$, but we are unable to
be exact.

\section{Square Quadruples}
We now consider the problem of finding $\{x_1,x_2,x_3,x_4\}$ with
all $6$ combinations $x_i^2 x_j^2+x_i^2+x_j^2=\square$ for $i \ne
j$. For ease of notation, we call the quadruple $\{x,y,z,w\}$.

Assume $x=e/f$, with $e$ and $f$ positive and having no common
factors. The rational solutions to $(x^2+1)t^2+x^2=\square$ can be
parameterised by $t=2xm/(x^2+1-m^2)$, so that
\begin{displaymath}
y=\frac{2 x m}{x^2-m^2+1} \hspace{2cm} z=\frac{2 x n}{x^2-n^2+1}
\end{displaymath}

Assuming that $m=g/h$, we find that $y^2 z^2+y^2+z^2=\square$
requires that
\begin{displaymath}
m^2 n^4+(m^4-4m^2+x^4+2x^2+1)n^2+m^2(x^4+2x^2+1)=\square
\end{displaymath}

Considering $m$ and $x$ as parameters, this quartic in $n$ has an
obvious solution when $n=0$, so is equivalent to an elliptic
curve. We find the following curve with integer coefficients,
$J^2=K^3+A K^2 + B K$ where
\begin{displaymath}
A=-2(e^4 h^4+2e^2 f^2 h^4+f^4(g^4-4g^2h^2+h^4))
\end{displaymath}
and
\begin{displaymath}
B=e^8 h^8 + 4e^6f^2h^8 - 2e^4f^4h^4(g^4+4g^2h^2-3h^4) -
\end{displaymath}
\begin{displaymath}
4e^2f^6h^4(g^4+4g^2h^2-h^4) +
f^8(g^8-8g^6h^2+14g^4h^4-8g^2h^6+h^8)
\end{displaymath}
with $n=J/(2f^2h^2K m)$.

We search this curve for integer points and generate $n$ from
these, giving the Eulerian triple $\{x^2,y^2,z^2\}$ .

We now need to find $w$ with $(x^2+1)w^2+x^2=\square$ assuming
that $x$ is a known fixed value. From before, $w=2 r
x/(x^2+1-r^2)$, which we substitute into
$(y^2+1)w^2+w^2=\square$, to find that $r$ must satisfy
\begin{displaymath}
y^2r^4+(2x^2(y^2+2)-2y^2)r^2+x^4y^2+2x^2y^2+y^2=\square
\end{displaymath}

This quartic has an obvious solution when $r=0$, so it equivalent
to an elliptic curve. We find the curve
\begin{displaymath}
V^2=U^3+(x^2(2y^2+1)+y^2)U^2+(x^4y^2(y^2+1)+x^2y^2(y^2+1))U
\end{displaymath}
with $r=V/(y(x^2(y^2+1)+U))$.

This elliptic curve can be easily transformed to have integer
coefficients and then searched for integer points, some of which
lead to values of $w$ different from $x,y,z$. The final test is
to compute $w^2z^2+w^2+z^2$ and to test whether this is square.

This methodology was coded using the simple multiple precision
UBASIC system, and run for several hours on a PC. The code finds
hundreds of examples with $5$ of the $6$ identities equal to a
square, but only one example with all $6$ square. This solution is
\begin{displaymath}
\left\{ 18^2 , \left(\frac{3}{5}\right)^2 ,
\left(\frac{8}{5}\right)^2 , \left(\frac{224}{107}\right)^2
\right\}
\end{displaymath}

\end{document}